
\magnification=1200

\def\proof{\smallskip\noindent{\bf Proof. }}
\def\remark{\medskip\noindent{\bf Remark}}
\def\example{\noindent{\bf Example }}
\def\definition {\noindent{\bf Definition}}

\def\O{{\cal O}}

\def\E{{\cal E}}

\def\P{{\bf P}}
\def\Q{{\bf Q}}
\def\P{{\bf P}}
\def\Z{{\bf Z}}
\def\C{{\bf C}}

\def\f{\varphi}
\def\ra{\longrightarrow}
\def\iso{\simeq}

\def\har#1{\smash{\mathop{\hbox to .8 cm{\rightarrowfill}}
\limits^{\scriptstyle#1}_{}}}

\font\ninerm=cmr10 at 10truept
\font\ninebf=cmbx10 at 10truept
\font\nineit=cmti10 at 10truept
\font\bigrm=cmr12 at 12 pt
\font\bigbf=cmbx12 at 12 pt
\font\bigit=cmti12 at 12 pt

\def\smalltype{\let\rm=\ninerm \let\bf=\ninebf
\let\it=\nineit  \baselineskip=9.5pt minus .75pt
\rm}
\def\bigtype{\let\rm=\bigrm \let\bf=\bigbf\let\it=\bigit
\baselineskip=12 pt minus 1 pt \rm}
 
{\bigtype
{\bf \centerline{Actions of Linear Algebraic Groups on Projective Manifolds}} 
{\bf \centerline{and Minimal Model Program }}}
\bigskip
\centerline{\it Marco Andreatta}

\bigskip
Dipartimento di Matematica, Universit\'a di Trento,

38050 Povo (TN), Italy

e-mail :  {\tt andreatt@science.unitn.it}

\null \noindent

\vskip 1 cm

Let $G$ be a connected linear algebraic complex group
which acts regularly and non trivially on a
smooth connected projective complex variety $X$ of dimension $n$.

In this paper we consider the following question: how does 
the $G$-action influence or even determine the structure of $X$?

As it is stand this is a too general question, thus we will
soon add some suitable assumptions; however even in this generality
we notice that $X$ is not minimal in the sense of Minimal Model Program (MMP).
In particular $X$ admits an extremal ray and an associated extremal 
(or Fano-Mori) contraction, $\f :X \ra Z$, which turns out to be
$G$-equivariant. 
It is thus natural to use the tools developed by the 
MMP, and the good properties
of the map $\f$, to get a classification of such varieties $X$.
This idea was first developed in an important paper by Mukai and Umemura 
(see [M-U]), where they studied smooth projective 3-folds on which 
$G= SL(2)$ acts with a dense orbit.
(A complete classification of such quasihomogeneous 3-folds
is contained in a paper of Nakano (see [Na2]); we refer the reader
also to a more recent work by S. Kebekus where the case
of singular 3-folds is also considered (see [Ke])).

Note also that if $X$ is actually homogeneous with the respect to $G$,
then $X$ is a Fano manifold and $X$ can be classified in term of Dynkin
diagrams. Fano manifolds are basic blocks of the MMP and 
moreover in this case there is a beutifull interplay between the 
representation theory of $G$ and the projective or differential geometry of $X$;
one aspect of it is the so called Borel-Weil-Bott theory. 
  
We want to propose a way to attack the general problem; however, to our knowledge, this way
works effectively only in the case when $G$ is a simple group, i.e. the simply connected
Lie group associated to a simple Lie algebra. In this case one can in fact perform 
many computations which seems hard or meaningless otherwise (for instance 
find the minimal non trivial irreducible representation).  

Thus we will also assume that $G$ is a (simply connected) simple Lie group and
we will define $r_G$ to be the minimum of the dimension of the 
homogeneous variety of the group $G$. That is $r_G$ is minimum codimension of
the maximal parabolic subgroup of $G$ (i.e. parabolic subgroup corresponding
to one node of the Dynkin diagram). It is easy to calculate $r_G$ if $G$ is simple
and this is done in section 1 (for instance $r_{SL(m)}= m-1$).

Then we first prove that if $r_G > n$ there is no such an $X$, that is the only possible
regular action is the trivial one, while if $r_G = n$ then $X$ is homogeneous. 
For instance if $n=3$ this says in particular that the only 
classical group acting non trivially on a smooth 3-fold are 
$SL(2),SL(3),SL(4), Sp(4)\iso Spin(5), SO(4)$ and in the last 3 cases $X$ is homogeneous; 
this special case was first proved in a paper of T. Nakano
(see [Na1]) which influenced the set up of this paper.

Then we classify all $X$ in the case $r_G = n-1$, see the theorem (4.1),
via the MMP. 
The special case in which $G = SL(n)$ was obtained first by T. Mabuchi
but in a complete different way. Namely he started with the classification of
$n$-codimensional closed subgroups of $SL(n)$, which follows from Dynkin's work, 
and he discussed the possible completions of their quotients.

Finally we begin to consider the case $r_G = n-2$; this is much more difficult and it seems
reasonable to make an additional general assumption. Namely we will assume that 
$X$ has an open dense orbit; such an $X$ is called a quasi homogeneous manifold.
As remarked above the case with $n= 3$ and $G = SL(2)$ was studied in [M-U]
and [Na2] while the case with $n=4$ and $G= SL(3)$ was recently settled by
Nakano, in [Na3],
with the method of computing the closed subgroup of codimension 
$4$ in $SL(3)$.
In the present paper we sketch an alternative proof of this
classification which makes use of the MMP.
Our hope is to extend in the future this method to the general case of $SL(n-1)$ and to the other
classical groups.

I was very much inspired by the papers of Mukai-Umemura, of Mabuchi and of Nakano 
which are quoted in the references.

This note was initiated during my visit at the University of Utah in the fall 
of 1997.  J. Koll\'ar suggested me to investigate in this direction using the 
Minimal Model Theory and in particular the new results I 
had recently with J. A. Wisniewski. The proof of the propositions (2.1)
and (2.2) were worked out with his help; I like to thank him for all this.

I also thank E. Ballico and P. Moeseneder for helpful discussions on this topic and the
organizers of the {\sl Hirzebruch 70} conference in Warsaw 
where this material was presented.

\beginsection 1. Definitions and preliminaries.

In this paper $X$ will always denote a smooth connected projective variety of dimension $n$.
We use the standard notation from algebraic geometry; more precisely
for the Minimal Model Program our notation is compatible
with that of [K-M-M]
while for the Group action and representation theory
it is compatible with that of [F-H].

\bigskip
\definition (1.1). A smooth projective variety $X$ 
will be said minimal (in the sense of the MMP) if $K_X$ is nef.

\smallskip
\proclaim Theorem (Mori-Kawamata-Shokurov) (1.2). 
Let $X$ be a smooth variety which is not minimal. Then there exists a map
$\f:X \ra Z$ into a normal projective variety $Z$ with connected fibers such that
$-K_X$ is $\f$-ample and $\f$ contracts the set of curves numerically equivalent
to a (rational) curve in a non trivial fiber. 

\smallskip
\definition (1.3). The map $\f:X \ra Z$ given in the above theorem
is called an extremal contraction or a Fano-Mori contraction.

\bigskip
Let $G$ be a connected linear algebraic complex group
acting regularly and non trivially on $X$. 
If $\pi :\tilde G \ra G$ is the universal covering map of $G$
then $\tilde G$ acts regularly and non trivially on $X$ through $\pi$. 
Hence we may and shall assume that the acting simple group 
is simply connected without loosing generality.

\proclaim Lemma (1.4). Let $X$ and $G$ be as above. 
Then $X$ is uniruled and in particular it is not minimal.

\proof On the generic point the action is not trivial, hence it is contained in an orbit 
which is unirational since $G$ is rational. Thus the generic point is contained in a 
rational curve $\subset X$. Therefore $X$ is uniruled and not minimal (for this
last statement see for instance [Ko], p. 189; the smoothness of $X$ is at this
point a necessary condition).

\smallskip
\definition (1.5).  If the action of $G$ is transitive on 
$X$ then $X$ is called a {\sl homogeneous manifolds}.
If $X$ has a dense open orbit then it is called a {\sl 
quasi-homogeneous manifold}.

\remark (1.6).If $X$ is homogeneous then $T_X$ is generated by
global sections and $-K_X$ is ample
(see for instance [Ko], (v.1.4)); in particular $X$ is a 
rational Fano manifold. 
If $X$ is quasi-homogeneous then $-K_X$ is effective;
this follows easily taking $n$ elements of the Lie algebra 
$Lie(G)$ such that their associated vector fields are linearly independent
at a generic point of $X$. The wedge product of these vector
fields gives a non trivial holomorphic section of $-K_X$.

\proclaim Lemma (1.7). Let $X$ and $G$ be as above.
Then there exists on $X$ a Fano-Mori contraction,
$\f:X \ra Z$, which is $G$-equivariant and $G$ operates on $Z$.

\proof The existence of $\f$ follows from the lemma (1.4)
and the above Mori-Kawamata-Shokurov theorem (1.2).
The equivariance of $\f$ follows from the following two facts:
on one end two curves which are carried one to another by the action of $G$ 
are numerically equivalent, on the other end $\f$ contracts 
all and only the set of curves 
in a ray, i.e. a set of curve all numerically equivalent
to a (rational) curve in a non trivial fiber.
Therefore take two points in a fiber and a curve passing through
these two points; this curve will be carried into another curve by the action of 
$G$ which is numerically equivalent to the first one and therefore it is contained in 
a fiber.

\bigskip 
\definition (1.8). Let us fix a simple, simply connected and connected Lie group $G$ 
and consider the set of all homogeneous manifolds (of dimension $>0$) with respect to this group. 
They are in a direct correspondence with the parabolic subgroups of $G$
(the isotropy subgroup in one point) which are in turn in direct correspondence
with the subsets of the nodes of the Dynkin diagram associated to the group $G$. 
We define $r= r_G$ to be the minimal of the dimensions of
the manifolds in this set or equivalently the {\sl minimal codimension of 
parabolic subgroups of $G$}.
A homogeneous variety which attains this minium will be called 
a {\sl minimal homogeneous variety} for the action of $G$. 
The minimal codimension will be attained at a maximal
parabolic subgroup, i.e. the ones corresponding to a single
node of the Dynkin diagram.

\smallskip
\example (1.8.1).
It is easy to check that 
if $G = SL(m)$ or $Sp(2s)=Sp(m)$ and $s \geq 3$ then $r_G=m-1$. 
If $G =SL(m)$ the parabolic subgroup $P$ is the
one corresponding to the first (or the last) node of 
the Dynkin diagram $A_m$;
if $G = Sp(2s)$ then $P$ is the
one corresponding to the first node of the Dynkin diagram $C_{s}$.
In both cases $G/P = \P^{(n-1)}$ and the actions comes from the identification 
$\P^{(m-1)}=\P(V)$, where $V$ is the m-dimensional irreducible representation
(this is called the standard homogeneous action).

Also if $G = Sp(4)$ then $r_G=3$ but in this case we have 
two different parabolic group of codimension $3$ which are the subgroup 
$P_1$ corresponding to the first node and $P_2 $ corresponding to the second one in the Dynkin
diagram; in this case $Sp(4)/P_1 = \Q^3$ and $Sp(4)/P_2 = \P^3$.

Note that $Spin(5) \iso Sp(4)$ and  $Spin(6) \iso SL(4)$; thus
when we consider the group $G = Spin(m)$ we will always assume that
$m\geq 7$ 

If $G = Spin(m)$ and $m\geq 7$ then $r_G = m-2$. If $m\not= 8$ the parabolic
subgroup $P$ is the one corresponding to the first node of the Dynkin diagrams 
$B_{(m-1/2)}$ or $D_{(m/2)}$, depending if $m$ is odd or even 
and $G/P \iso \Q^{(m-2)} \subset \P^{(m-1)}$.
If $G= Spin (8)$ in principle we will have two minimal homogeneous varieties
(spinor varieties) of dimension 6
(corresponding to each of the two last nodes)
but they are both isomorphic to $\Q^6$. 

Similarly one can easily compute $r_G$ when 
$G$ is an exceptional group. 
For instance if $G = G_2$ then $r_G = 5$ and $G/P = \Q ^5 \subset \P^6$.

\medskip
\definition (1.9). Let $X = G/P$ be an homogeneous variety where $G$ is a 
simply connected simple group and $P$ is a parabolic subgroup.
A vector bundle $E \ra X=G/P$ is called $G$-homogeneous or simply 
homogeneous if there exists an action of $G$ over $E$ such that the following diagram commutes

\centerline{$G\times E \ \ \ra \ \ \ E$}
\centerline{$\downarrow  \ \ \ \ \ \ \ \ \ \ \ \ \ \ \ \ \  \downarrow$}
\centerline{$G\times (G/P) \ra G/P$}

\remark (1.10). It is evident from the definition that the tangent bundle
of $X$ is homogeneous. 

One can prove that a vector bundle $E$ on $X =G/P$ is homogeneous 
if and only if one of the following conditions
holds:
\item{i)} $\theta_g^*E \iso E$ for every $g \in G$ ; $\theta_g$ is the automorphim
of $X$ given by $g$.
\item{ii)} There exists a representation $\rho : P \ra GL(r)$ such that
$E \iso E_{\rho}$, where $E_{\rho}$ is the vector bundle with fiber $\C^r$ coming
from the principal bundle $G \ra G/P$ via $\rho$.

\beginsection 2. Points which are fixed by the action of $G$.

In the first part of this section we enlarge slightly our set up:
namely we will have an action of a connected linear algebraic group $G$ on 
a variety $Z$ with normal singularities. 
We want to discuss how the existence of a fixed point by the action of 
$G$ determines the structure of $Z$, at least in a Zariski open subset.
The next proposition (2.1) is what is needed for our pourposes; 
after proving it we actually found that the following very nice and general result of 
Ahiezer implies easily our proposition.

\proclaim Theorem (2.0). ([Ah] theorem 3, [H-O] for the analytic case)
Suppose that a connected linear
algebraic group $\tilde G$ acts effectively on a complete normal variety $X$ with an open orbit
$\Omega$ such that $A := X \setminus \Omega$ contains an isolated point $p$.
Then the variety
is a projective cone over a homogeneous variety $P$ with the respect to the maximal 
reductive subgroup of $\tilde G$. Moreover $E$ is either $p$ or
$p \cup P$.  

We think however that it can be useful to provide here both a direct and
more geometric proof of the proposition as well as an argument showing that it is a consequence of
the above theorem. The main steps of both proofs grew from suggestions by J\'anos Koll\'ar.

\medskip
We first recall some standard notations and definitions
which will be used in the proof.
Let $Z$ be a variety with normal singularities and $z\in Z$ be point on it; 
with $gr(\O_{Z,z})$ we denote as usually the graded ring associated to the 
maximal ideal $m_z$, 
that is $:=\bigoplus_k m_z^k/m_z^{k+1}$. Then $V_z:= Spec [gr(\O_{Z,z})]$ is the tangent cone
to $Z$ at $z$ and $P_z:= Proj [gr(\O_{Z,z})]$ is the projectivized tangent cone, 
or equivalently the exceptional divisor of the blow-up of $Z$ at the point $z$.
The surjective morphism 
$\bigoplus Sym^k(m_z/m_z^2) \ra \bigoplus m_z^k/m_z^{k+1}$
gives the embedding of $P_z$ into the projectivization 
of the Zariski tangent space in $z$, which is equal to $\P^N$ where $N+1= dim (m_z/m_z^2)$.
This embedding is linearly normal and it is given by the global sections
of the twisting line bundle $L$ on $P_z$.

Assume that a simple group $G$ acts on $Z$ and $z$ is a fixed point for this action;
thus $G$ acts on the $\C$-vector spaces $m_z^k/m_z^{k+1}$
for every $k$ and therefore $G$ acts on $P_z$.
Assume also that $P_z$ is an homogeneus space under this action; 
that is $P_z = G/ P$ with $P$ a parabolic
subgroup of $G$. Then $L$ is a homogeneus line bundle
and it is determined by a character $\lambda : P \ra \C^*$; 
moreover, since $L$ is spanned, the character $\lambda$ is a positive multiple
of a highest weight.
By the Borel-Weil-Bott theorem for positive integer $k$ 
the vector space $H^i(P_z, L^k)$ are zero for $i \not= 0$ and
$V_{\lambda ^. k}: = H^0(P_z, L^k)$ is an irreducible $G$-module; i.e.
it is the irreducible reppresentation with highest weight $\lambda ^. k$. 
In particular, since we have the $G$-equivariant inclusion 
$m_z^k/m_z^{k+1} \subset H^0(P_z, L^k)$, these two spaces are the same one
(and the above described embedding $P_z \ra P^N$ is projectively normal).

The tangent cone $V_z$ has a natural G-structure which is induced by the action
on $P_z$ (and on $L$); with this action $V_z$ has two orbits, namely the vertex,
which is a fixed point, and its complement which is an open orbit.

\proclaim Proposition (of the cone) (2.1). Let $G$ be a simple Lie group acting on a normal 
projective variety $Z$. Assume that $z$ is a fixed point and that 
$P_z$ is a homogeneus variety with the respect to the action of $G$.
Then $Z$ has a Zariski open set which is $G$-equivariantly isomorphic
to the (affine) cone over $z$.

\proof 
Let for brevity $R = \O_{Z,z}$ and $m = m_z$; let also $k$ be a positive integer.
By the complete reducibility of the group $G$ there is a unique complement of $m^{k+1}$
in $m^k$; let us call it $V_k$. Note that $V_k \subset m^k \subset R$; thus
if we define $A: =\bigoplus_k V_k$ we have that $A \subset R$.

We will first prove that $A$ has ring structure (and thus that it is equal to $gr(\O_{Z,z})$).
This follows if we prove that the image of the natural map
$$\f _k : S^k (V_1) \ra m^k$$
is exactly $V_k$, and thus that $A$ is a quotient of a polynomial ring.

$Im (\f_k):=W_k$ is contained in $m^k$ and it is actually mapped injectively into
$m^k/ m^{k+s}$ for sufficiently large $s$. But on one side $W_k$ is a representation
of $G$ with highest weight $\leq k$, since $S^k (V_1)$ is a representation 
with highest weight $k$. On the other $m^k/ m^{k+s}= V_k + V_{k+1} +....V_{k+s-1}$,
the direct sum of irreducible distinguished representation $V_i$ of highest weight 
$i$, for $i = k,..., k+s-1$. This implies that $W_k= V_k$.

\smallskip
The ring inclusion $A \subset R$ gives a G-invariant map $\rho: Spec(R) \ra Spec(A)$;
we will now prove that it is an isomorphism. 

Note that $\rho$ is unramified at $z$: in fact this is equivalent to say that 
the sheaf of relative differential
$\Omega_{Spec(R)/Spec(A)}$ is zero at $z$. But the stalk of the sheaf in $z$ is equal to
$\Omega_{R/A}$, which fits in the exact sequence
$$\Omega_{A/\C} \ra \Omega_{R/\C} \ra \Omega_{R/A} \ra 0$$
and the first map is the isomorphim $m_A/m_A^2 \ra m/m^2$
(where $m_A$ is the maximal ideal of $A$). 

Note also that, since $G$ acts transitively on 
$Spec(A)\setminus v$, we have that 
$\rho _{|Spec(R)\setminus z}: Spec(R) \setminus z \ra Spec(A)\setminus v$ 
is an \'etale cover.
In fact, by dimension reasons, it is generically finite and it is also
generically flat. Using the transitivity of the action
we conclude then that it is finite and flat and that it has no ramification points.

To describe this cover we need the following lemma of topology; it is quite straightforward
but we will give a proof for the reader's convenience.   

\proclaim Lemma (2.1.1). Let $P$ be a smooth complex compact variety which is
simply connected;
let $L$ be a line bundle on $P$. Denote by $L^*$ the total space of the bundle $L$ 
minus the zero section; $L^*$ is a $\C^*$ bundle over $P$.
Let $k$ be the unique element in $H^2(P,\Z)$ which is not divisible
(i.e. is not a multiple of any element in $H^2(P,\Z)$)
and such that $rk = c_1(L)$.
Then $\pi _1(L^*) = \Z /r\Z$ and if $K$ is a line bundle which has Chern
class $k$ then $K^*$ is the universal cover of $L^*$ (the covering
map is described in the course of the proof). 

\proof
We start with the exact sequence of homotopy

$$0 \ra \pi_2(L^*) \ra \pi_2(P) \ra\pi_1(\C^*) =\Z \ra \pi_1(L^*) \ra \pi_1(P) =1.$$

Since $\pi_1(P)=1$ then, by Hurewicz's theorem, we have $H_2(P,\Z)=\pi_2(P)$
and, by the universal coefficient theorem for cohomology, $H^2(P,\Z)$
has no torsion.
The above map $\pi_2(P)= H_2(P, \Z) \ra \Z$ is nothing else 
then the first Chern class $c_1(L)\in H^2(P,\Z)$.

So $\pi_1(L^*)$ is a co-torsion of the map $c_1(L)$ or, in other
words, its a cyclic group of rank $r$, where $r$ is the integer in
the statement. 
(This is a standard fact about
abelian groups with no torsion: i.e. if we have an exact sequence
$ A  ^{g \atop \ra}  \Z \rightarrow B \rightarrow 1$, where $g\in  A^*$ is equal to $rc$ with
$c$ not divisible, then $B = \Z/ r\Z$). 

Let $K$ be the line bundle on $P$ which has Chern
class $k$; it exists since $rk$ is a Chern class and $Pic(P) \ra H^2(P,\Z)$ 
has a torsion free cokernel.

Now let $K'$ be the pullback of $K$ on 
the total space of $L$. Then the zero section of $L$, $P_0$, is 
a section of $rK'$ and thus we can construct a cyclic covering of $L$
branched along $P_0$ with degree $r$. Then we
throw away the branch locus $P_0$ and we obtain a covering 
which is the universal one, by the previous consideration.

\smallskip
The lemma concludes the proof of the proposition: in fact it describes 
the etale cover $\rho_{SpecR \setminus z}$ and it implies in particular that 
if $\rho$ is not an isomorphism then it is ramified at $z$, which contradicts
what we have proved above.

\bigskip
We now show how to prove proposition (2.1)
(actually a slightly stronger version of it) using the theorem of Ahiezer.

\proclaim Proposition (2.2). Suppose that a connected reductive linear
algebraic group $G$ acts effectively on a complete normal variety $Z$.
Then the followings are equivalent:
\item{a)} There exists a fixed point $z$ such that 
$P_z$ is a $G$-homogeneus variety.
\item{b)} $Z$ has an open orbit $\Omega$ 
and $A := Z\setminus \Omega$ contains an isolated point $z$.
\item{c)} $Z$ is a projective cone over a homogeneous variety $P$ with the respect to $G$.

\proof It remains to prove that a) implies b), the rest following from (2.0).
In the assumption of a) let thus $O$ be the smallest orbit, with the respect to
the dimension, in the closure of which $z$ is contained. 
The tangent cone of this orbit and of its closure is $G$ invariant.
But, since by assumption $P_z$ is $G$-homogeneous, this implies that the orbit has dimension
equal to the dimension of $Z$.

\bigskip
A first application of the above proposition will give the next result.

\proclaim Lemma (2.3). Let $X$ be a smooth projective variety and $G$ 
a simple, simply connected, connected linear group acting non trivially on $X$;
let $r_G$ be the integer defined in (1.8) and $n=dimX$. 
If $n \leq r_G$ there are no 
fixed points on $X$. 
If $r_G = (n-1)$ then $X$ has no fixed points unless 
$G = SL(n)$ or $Sp(n=2s)$, $X = \P^n$ and the action is the one which extends
the standard $SL(n)$ or $Sp(n)$ action on $\C^n$ via the inclusion
$\C^n \ra \P^n, (z_1,...,z_n) \ra (1,z_1,...,z_n)$
(equivalently the action is induced from the homomorphism 
$g \ra \left(\matrix{1&0\cr 
0&g\cr}\right)$ from $SL(n)$ or $Sp(n)$ to 
$PGL(n+1)$).

\proof Note first that the dimensions of the irreducible representations
of $G$ are strictly bigger then $ r_G$: in fact for every irreducible 
representation $V$ there is a unique closed orbit in $\P(V)$ which
is the homogeneous variety corresponding to the parabolic subgroup
perpendicular to the weight of the representation. 
Moreover if the minimal of such dimension is equal to $(r_G+1)$ then $G = SL(m)$ or 
$Sp(m)$ and $V$ is the standard representation; in this case
the action on $\P(V) = \P^{m-1}$ is homogeneous.

Assume that $r_G \geq n$
and that $x\in X$ is a fixed point; then $G$ acts on the 
tangent space $T= T_x X$ and by the above observation this has to be the trivial 
representation. Let $m_x$ be the maximal ideal of $\O = \O_x$, the local ring
of germs of regular functions near $x$; then $G$ acts trivially on $m/m^2 =T^*$
and on $m^k /m^{(k+1)} = S^k (m/m^2)$.
Using inductively the exact sequences 
$$0 \ra m^k/m^{(k+1)} \ra \O /m^{(k+1)} \ra \O /m^k \ra 0$$
and the fact that $G$ is a reductive group we have that $G$
acts on $\O/m^k$ trivially for all $k>0$.
Thus $G$ acts trivially on the completion ${\hat {\O}}$, hence trivially
on $\O$. This implies that $G$ acts trivially on $X$ itself.

After noticing that $G$ acts trivially on 
$T$, one can conclude alternatively via the Luna's etale slice theorem
as in the next lemma (2.4).

Assume now that $r = (n-1)$ and let $x \in X$ be a fixed point. 
If $G= Spin(m)$ ($m \geq 7$) or an exceptional group then the above proof applies; 
i.e. the action of $G$ on the tangent space at $x$ must be trivial.
In the other cases we can apply the proposition (2.1) 
(or the proposition (2.2)) since
the action of $G$ on $P_x:= Proj [gr(\O_{X,x})] = \P^{(n-1)}$
is transitive. Thus $X$ contains a Zariski open set which is
isomorphic to the cone over a smooth point.
Let $Y = \P(O_{\P^{(n-1)}}(1) \oplus \O_{\P^{(n-1)}})$
with the action of $SL(n)$, or of $Sp(n)$, determined by the homogeneos
action on $\P^{(n-1)}$ and on the homogeneous bundle
$\O(1)$. Since $X$ and $Y$ have the same open 
orbit there is a birational map $f : X ---> Y$ induced by identifying the orbit.
By the Hironaka's theorem we can resolve the indeterminacy
locus $J$ of $f$ by succesive blowing-up along smooth centers, 
$\sigma :X' \ra X$.
Since $I$ is a $G$-stable closed subset of codimension greater or equal then $2$,
$\sigma$ is a composition of blow-ups at points.
Let $\tau = f \circ \sigma$; since the indeterminacy locus
of $\tau ^{-1}$ is also a $G$-stable closed subset of codimension $\geq 2$
and since there are no such subset on $Y$ we have that 
$\tau$ is an isomorphism and that $X$ is the blow-down
of $Y$ along a divisor to a point, i.e. $X = \P^n$.

\proclaim Lemma (2.4). Let $G = SL(n-1)$ acting with a dense open orbit on a 
n-fold $X$. Then there is no fixed points.

\proof If $n=2$ this is the lemma (1.2.2) in [M-U]. Therefore we assume that $n\geq 3$
and that, by contradiction, $x$ is a fixed point. Then we have an induced linear action
of $G$ on $T_{X,x}$, i.e. an n-dimensional representation of $G$.
These are of three types, namely if $A \in SL(n-1)$
$$A \longrightarrow (A,1),\hbox{ or } (^t A^{-1},1),\hbox  {or }  I;$$
in particular there are no $n$-dimensional orbit on $T_{X,x}$ in any of these 
there cases.

On the other hand we can apply the Luna's etale slice theorem (see [Lu]); this says  that
there exists a $G$-stable affine subvariety $V$ containing $x$ and an etale $G$-equivariant
morfism $V \ra T_{X,x}$. This is a contradiction since, by assumption, $X$
has a $n$-dimensional orbit.

\smallskip
Actually the following more general result holds; it was proved for
$n=3$ in [Na], here we adapt this proof (or the one of 1.2.2 in [M-U])
to the general case.

\proclaim Lemma (2.5).  Let $G$ be any reductive 
group acting with a dense open orbit on a 
projective variety $Z$ and assume that $x$ is a fixed point. Then
$m_x/m_x^2$ does not have nonzero invariants.

\proof Assume by contradiction that there exists a non-zero invariants
$f\in m_x/m_x^2$. Let $U = Spec(A)$ be a $G$-invariant affine neighborhood
of $x$ (this exists since $G$ is geometrically reductive). Let $\bar f$ be
a lifting of $f$, i.e. $\bar f \in \Gamma(U, \O_U)$ is  such that $\tau(\bar f) = f$
where $\tau :\Gamma(U,\O_U) \ra \O_y \ra \O_y/m_y$. Let $V$ be a finite dimensional
$G$-invariant vector subspace of $A$ containing $f$; this exists by Borel [Bo]
(it can be defined as the vector subspace of $A$ generated by $\{ g\circ \bar f| g \in G\}$
which is of finite dimension). Since $\bar f(x) = 0$ we have that $\tau (V)  \subset m_x/m_x^2$. 
The image $\tau (V)$ contains a non zero $G$ invariant hence $V$ contains a 
$G$ invariant. Since $V$ and $ m_x/m_x^2$ are finite dimensional,
and $G$ is lineraly reductive, the image $\tau (V)$ is a direct summand of $V$;
hence $V$, in particular $A$, contains a non zero $G$ invariant $F$. Since $G$
has an open orbit the invariant $F$ should be constant. Since its value on $x$ is zero 
it is constantly zero which is a contradiction.

\beginsection 3. A starting point.
 
Our main goal will be a classification 
of smooth connected projective varieties with a non trivial action of 
a simple group $G$ which has the number $r_G$ ''big enough''
with the respect to the dimension of $X$. 
The following easy result seems to be a good starting point.

\proclaim Proposition (3.1). Let $G$ be a simple, simply connected and connected Lie
group acting on a smooth projective variety $X$.
If $r_G > n$, then the action of $G$ on $X$ is the trivial one.
If $r_G = n$ and $G$ is acting non trivially then $X$ is homogeneus.
In particular if $G = SL(m)$ or $Sp(m)$ acts on a
smooth connected projective varieties $X$ of dimension $n < m-1$
then this action is trivial;
if $n = m-1$ then $X = \P^{(m-1)}$ and the action is the transitive standard one
apart for the case $G = Sp(4)$ where we have both $\P^3$ and $\Q^3$ as 
homogeneous variety of dimension $3$.
If $G= Spin(m)$ with $m \geq 7$ acts on a smooth connected projective varieties 
$X$ of dimension $n < m-2$ then this action is trivial; 
if $n=m-2$ then $X = \Q^{(m-2)}$ and the action is the transitive standard one.   

\proof
The proof of the proposition follows easily from the assumption on $r_G$
and the fact that on one side 
there are closed orbits on $X$ (see for instance [Bo], 1.8) and on the other
$X$ has no fixed point as showed in the first part of the lemma (2.3).

\remark (3.1.1). The special case $n=3$ of the proposition gives the main theorem of [Na1]; 
the proof of the first part of the lemma (2.3) above, the one not using the Luna's theorem, 
is the same as the one of lemma 5 in [Na1].

\beginsection 4. Minimal Model Program on manifolds with a G-action.

Let $X$ be a smooth projective manifold of dimension $n$ 
and $G$ a simple, simply connected and connected Lie group acting non trivially on $X$. 

We use the notation and the approach of the previous section,
passing to the next step; namely we assume that $r_G = (n-1)$.

More precisely we will prove the following theorem, the first part of
which was proved in [Ma3] with different methods.

\smallskip 
\noindent{
{\bf Theorem (4.1)} {\sl If $G= SL(n)$ then $X$ is isomorphic
to one of the following varieties; the action of 
$G$ is unique for each case and it is described in the course of the proof (see also [Ma]):
\item{1)} the complex projective space $\P^n$,
\item{2)} $\P^{(n-1)} \times R$, where $R$ is a smooth projective curve,
\item{3)} The projective bundles $\P(\O_{\P^{(n-1)}}(m) \oplus \O_{\P^{(n-1)}})$ with $m>0$, 
\item{4)} if $n=2$ we have moreover a different action from the above one on 
$\P^1 \times \P^1$ and $\P^2$,
\item{5)} if $n= 3$ we have moreover the projective bundle $\P(T(\P^2))$,
\item{6)} if $n=4$ we have moreover the smooth 4-dimensional quadric which is isomorphic
to $Gr(2,4)$, the Grassmannian of 2-planes in $\C^4$.

\smallskip  \noindent
If $G= Sp(n)$ then $X$ is isomorphic to one of the following
varieties and the action of $G$ is unique for each case.
\item{1)} the complex projective space $\P^n$,
\item{2)} $\P^{(n-1)} \times R$, where $R$ is a smooth projective curve,
\item{3)} The projective bundle $\P(\O_{\P^{(n-1)}}(m) \oplus \O_{\P^{(n-1)}})$ with $m>0$, 
\item{4)} if $n= 4$ we have moreover $\Q^4$, 
the homogeneous variety which is the quotient of $Sp(4)$ by the Borel subgroup
(which has two structure of a $\P^1$-bundle over $\P^3$ and over $\Q^3$), 
$\Q^3 \times R$, where $R$ is a smooth projective curve and the projective bundles
$\P(\O_{\Q^3}(m) \oplus \O_{\Q^3})$ with $m>0$.

\smallskip \noindent
If $G = Spin(n+1)$ with $n \geq 6$ then $X$ is isomorphic to one of the following
varieties and the action of $G$ is unique for each case.
\item{1)} the complex projective space $\P^n$,
\item{2)} the complex projective quadric $\Q^n \subset \P^{(n+1)}$,
\item{3)} $\Q^{(n-1)} \times R$, where $R$ is a smooth projective curve,
\item{4)} The projective bundle $\P(\O_{\Q^{(n-1)}}(m) \oplus \O_{\Q^{(n-1)}})$ with $m>0$.
}}

\proof The proof of the theorem will be reached in a number of steps
which are similar for all the three groups. 

\proclaim Lemma (4.2). Let $X$ and $Y$ two manifolds on which a simple
group $G$ acts in the hypothesis of the theorem (i.e. =$r_G+1 =dimX =dimY$).
Assume that $X$ and $Y$ have each a dense open orbit which are $G$ isomorphich,
then $X \iso Y$ unless $G = SL(n), Sp(n)$, $Y= \P^n$ and $X = \P(\O(1) \oplus \O)$.

\proof (See also the last part of the proof of (2.3)). 
Since both $X$ and $Y$ are completion of the same open dense
orbit there is a birational map $f : Y ---> X$ induced by identifying the orbit.
If $Y = \P^n$ let consider the blow-up of the fixed point 
$\sigma : Y' \ra Y$ and take instead of $f$ the composition $g = f \circ \sigma.$
This map is defined in codimension 1, since both $X$ and $Y$ has minimal closed orbits
of codimension $1$ and no fixed point, thus it is an isomorphism. 

\bigskip
Let us now run the Minimal Model Program to classify $X$; in the following $\rho(X)$
will denote the Picard number of $X$.

\medskip\noindent
1-st Step. Assume that $\rho(X) \geq 2$ and let $\f :X \ra Z$ be the contraction of 
an extremal ray. 

a) If $\f$ is birational then, by the $G$-equivariant property of $\f$
and our assumption on $r$, 
it must be divisorial and the divisor has to be contracted to a 
point. Moreover the exceptional divisor $E$ is isomorphic to $\P^{(n-1)}$,
respectively to $\Q^{(n-1)}$; here the two cases depends on whether 
$G= SL(n), Sp(n=2s)$ or if $G = Spin(n+1), n\geq 6$, unless $G= Sp(4) \iso Spin(5)$
in which both are possible.
Since it is an extremal contraction the conormal bundle
of the exceptional locus is $N^* = \O(k)$ with $1 \leq k \leq n-1$, respectively
$1 \leq k \leq n-2$.

We can thus apply the cone's proposition (2.1) (to $z\in Z$); 
this gives that $X$ is a completion of the open variety
$V(E, N^*) = Spec (\bigoplus_h \O(hk))$. Note that the 
open orbit is isomorphic to $G/K$ where $K$ is the kernel of the character
map $\rho : P \ra \C^*$ associated to the homogeneuos line bundle
$\O(k)$, $P$ is the parabolic subgroup associated to
$\P^{n-1}$, resp. $\Q^{n-1}$.

One possible completion is $X_k=\P(N^* \oplus \O)$ which has an open orbit 
isomorphic to $G/K$ and two closed orbit isomorphic to $\P^{(n-1)}$, 
respectively $\Q^{(n-1)}$.
But, by the above lemma (4.2), this is actually the only one except if 
$k=1$ and $G = SL(m)$ or $Sp(m)$, where $X_1$ can be actually blow-down to 
$\P^n$. In this case there are thus two possible completions
(actually $\rho (\P^n) = 1$ and thus it 
will appear in the proper place in the second step).

\smallskip
b) Let $\f$ be of fiber type and consider the induced action of $G$ on $Z$. 
By our assumption either this action is trivial or $Z = \P^{(n-1)}$ if $G= SL(n)$ or $Sp(n=2s)$, 
respectively $\Q^{(n-1)}$ if $G = Spin(n+1)$. 

In the first case,
since any fiber of $\f$ is an orbit, we must have that $dim Z = 1$ and $X = \P^{(n-1)}
\times Z$, respectively $\Q^{(n-1)} \times Z$, with the $G$-action factorizing 
to the product of the standard
homogeneous one on $\P^{(n-1)}$, respectively on $\Q^{(n-1)}$, and the trivial one on $Z$, 
except possibly for $n=2$. This follows for instance by the more general
theorem 1.2.1 in [Ma2]; for the reader convenience we outline his proof
in this case. Namely take a point
$p_0\in X$ and let $H$ be the isotropy group of $G$ at $p_0$.
Let $Z_1 = \{ p\in X : G_p = H\}$, where $G_p$ is the isotropy group of $G$ at $p$.
Then one can define a regular map $\tau : G/H \times Z_1 \ra X$ by 
$\tau (gH,p) = g^.p$. It is straightforward to see that this map is well defined,
injective and $G$-equivariant. Moreover, by the Zariski's main Theorem,
it is an algebraic $G$-equivariant isomorphism. This gives our claim after noticing
that $G/H \iso \P^{(n-1)}$,  respectively $\Q^{(n-1)}$, and that $Z_1 = X/G= Z$.

If $n = 2$ and $G = SL(2)$ then we have an other case which comes from 
the diagonal action of $SL(2)$ on $\P^1 \times \P^1$. 
It is straightforward to prove that
there are no other actions of $SL(2)$ on the smooth two dimensional quadric.

In the second one $\f$ is an equivariant $\P^1$-bundle over $\P^{(n-1)}$,
respectively $\Q^{(n-1)}$:
in fact the action on $Z$ is homogeneous and thus the fibers are 
all equidimensional and there are no reducible or double fibers.
Thus $X = \P(E)$ with $E$ a rank 2 vector bundle on $Z$; 
$E$ is homogeneous since the action is $\f$ equivariant. Therefore,
after normalizing if necessary, $E= \O(s) \oplus \O$ with $s\geq 0$ and, in the case
$n=3$ and $G= SL(3)$ one can also have $E = TZ$  while in the case
$n=4$ and $G= Sp(4)$ $E$ can also be the nullcorrelation bundle
on $\P^3$ or the spinor bundle on $\Q^3$. 

If $E= \O(s) \oplus \O$ we have a decomposition
of $X$ into three orbits. Two isomorphic to $\P^{(n-1)}$, respectively
$\Q^{(n-1)}$ (the sections at infinity
and the zero section) and a open dense orbit isomorphic to $G/S$ where
$S$ is the kernel of the character
map $\rho : P \ra \C^*$ associated to the homogeneuos line bundle
$\O(s)$, $P$ being the parabolic subgroup
associated to $\P^{n-1}$, resp. $\Q^{n-1}$. 
The fact that this is the unique action on $X$ can be proved as above
with the exception $n=2$ and $s = 0$  
(note that the section at infinity can be contracted so we can apply the cone's proposition).

If $n=3$ and $E = T\Z$ it is well known that $X = \P(T\P^2)$ is the homogeneous variety
$G/B$ where $B$ is a Borel subgroup of $SL(3)$ which corresponds in taking all
the Dynkin diagram $A_3$ (or equivalently the kernel
of the two dimensional representation of $H$ associated to the tangent bundle); 
it is the unique closed orbit of the adjoint representation of $SL(3)$.

If $n=4$ and $E$ is either the nullcorrelation bundle
on $\P^3$ or the spinor bundle on $\Q^3$ then $X= \P(E) = Sp(4)/B$ where $B$ is a Borel subgroup.

\medskip \noindent
2-nd Step. Assume finally that $\rho(X) =1$, i.e., since it 
has an extremal ray, $X$ is a Fano manifold. 

If $X$ is homogeneous then we can just look at the list of parabolic subgroups
of codimension $n$
corresponding to one node of the Dynkin diagram, in order to have  $\rho(X) =1$. 

If $G= SL(n)$ we have only one possibility for $n=4$, 
namely $X = SL(4)/Q$ where $Q$ is the parabolic subgroup corresponding to
the second node of the Dynkin diagram $A_4$. It is the unique orbit of the irreducible
representation of $SL(4)$ into $\Lambda ^2 \C^4$ and it is isomorphic to the Grassmanian
of planes in $\C^4$, i.e. the smooth $4$ dimensional quadric. 

If $G= Sp(n)$ or $Spin(n+1)$ with $n\geq 6$ there is no homogeneous manifold of dimension $n$
with  $\rho(X) =1$.

\smallskip
If $X$ is not homogeneous and has no fixed point then it must have a closed orbit $H$ 
whic will be isomorphic to $\P^{(n-1)}$, respectively $\Q^{(n-1)}$. 
Let $L$ be a positive generator of $Pic(X)$; then $H= mL$. Since $H$ is effective
$m> 0$; then it is well known 
that a smooth projective variety with an ample section isomorphic
to $\P^{(n-1)}$, respectively $\Q^{(n-1)}$, has to be isomorhic to 
$\P^n$ (if $n= 2$ we can have also $\P^1 \times \P^1$, this has however $\rho(X) =2$
and thus it was considered above) , respectively to $\P^n$ or to $\Q^n$.

So if $G = SL(n)$ or $Sp(n)$, the last with $n\not= 4$,
then $X$ has to be $\P^n$ and it contains
the closed orbit $H \iso \P^{(n-1)}$ as a linear subspace except for
$n= 2$ in which case the orbit can be a conic $\iso \P^1$. If the orbit is
linear then $X$ contains as an open Zariski subset the total space of the normal
bundle. Thus the action on this open subset is fixed (by the action on the orbit)
and as discussed above it is unique (see the lemma (4.2)). 

If $n=2$ we have another non trivial action : namely the induced action on $\P^2 = \P (\C^3)$
by the  $3$-dimensional irreducible representation
$\alpha _3 :SL(2,\C) \ra GL(3,\C)$. 
It is straightforward to prove that there are no other actions of $SL(2)$ on $\P^2$.

If $G = Sp(4)$ then we have the above case when $H \iso \P^3$ but we can have also
$H \iso \Q^3$. Then $X$ can be either $\P^4$ or $\Q^4$; the action is described in the following
if we think at $G$ as $Spin(5)$.

If $G = Spin(n+1)$ then we have an action on $X = \Q^n$
given by the embedding  $Spin(n+1) \ra Spin(n+2)$ and one can prove 
that this is the only 
possible action; there is a closed orbit, isomorphic to the 
$(n-1)$-dimensional quadric and a open orbit.
If $X = \P^n$
the action is coming from the canonical action of $G$ on
$\C^{(n+1)}$ and $X$ has two orbits: a closed one, isomorphic to the 
$(n-1)$-dimensional quadric and
a open one isomorphic to $X_2 = Spin(n+1) / S( O(1) \times O(n))$.

\remark (4.3). Results as in the theorem can be obtained 
if $G$ is an exceptional simple group and $r_G = (n-1)$.

\beginsection 5.  Fourfolds which are quasihomogeneous under the action of $SL(3)$.

The next step will be the case $r_G= (n-2)$, so for instance the case  
$G= SL(n-1)$ and $X$ of dimension $n$. 

If $n= 3$, $G = SL(2)$ and $X$ quasi homogeneous this was studied in a series of papers
starting with the one of Mukai-Umemura (see [M-U]); the complete results are collected in 
a paper of Nakano (see [Na2]). Here we will assume that $G= SL(3)$ and that 
$X$ is of dimension $4$ and quasi homogeneous with the respect to $G$. The classification
of such $X$ was already achieved by Nakano in [Na3]
by computing the closed subgroup of codimension $4$ in $SL(3)$.
Here we sketch an alternative proof, similar to the ones in [M-U] and [Na2]
for the 3-folds, with the aim to extend in the future
this method to the general case of $SL(n-1)$ and to the other
classical group; for this pourpose some results are given in general.
On the other hand many details, as an explicit description of the action or of
the closed subgroups of codimension $4$, are not worked out
and we refer the interested reader to the paper of Nakano.

\proclaim Lemma (5.1). There are no homogeneous n-fold with $Pic(X) =\Z$
for the action of $SL(n-1)$.

\proof $X$ should be the quotient of $SL(n-1)$ by a parabolic
subgroup corresponding to one node of the Dynkin diagram. There are no grassmanian
in $\C^{(n-1)}$ of dimension $n$.

\proclaim Lemma (5.2).  Let $G = SL(n-1)$ acting with a dense orbit on a n-fold $X$. 
Then, if $X$ is not homogeneous, there exists a closed orbit of dimension $(n-2)$
unless possibly if $n=4, 5$ and $X$ contains an orbit isomorphic to
$\P(T(\P^2))$, respectively $\Q^4$.
Moreover this orbit is isomorphic
to $\P^{(n-2)}$ and its normal bundle is homogeneous. 

\proof Since $X$ has no fixed point, see (2.4), it must have a closed orbit of dimension $\geq 1$;
since the group is $SL(n-1)$ the orbit has to be of dimension $\geq (n-2)$.
If it is of dimension $(n-2)$ then it is $\P^{(n-2)}$, if it has dimension
$(n-1)$ then it has to be one of the variety in the list of the Theorem in the previous section.
In particular the only homogeneous one are the one giving the exceptions in the lemma.

The fact that the normal bundle is homogeneous follows since it is a quotient
of the tangent of $X$ restricted to the orbit and the tangent of the orbit; these two bundles
are homogeneous on the orbit since the differential of every transformatrion in $SL(n-1)$ 
defines an automorhism of these bundles.

\bigskip
From now on we will assume that $X$ is a smooth $4$-fold quasihomogeneous
with the respect of a $SL(3)$ action and we will run the MMP on $X$.

We assume first that $\rho (X) \geq 2$ and let $\f:X \ra Z$ be the contraction of an 
extremal ray; if $\f$ is birational 
let also $E$ be the exceptional locus. As in the previous section, by
the $G$-equivariant property of $\f$ and the fact that $r_{SL(3)} =2$, we can, a priori, 
have only the following cases.

\smallskip \noindent
a)  $\f$ is birational $dim E = 3$ and $\f(E)=z$ is one point; in this case
$E$ is a $3$-dimensional del Pezzo variety with an $SL(3)$ action 
induced by the one of $X$. In particular $E$ cannot have fixed point.

First note that $E$ has to be smooth: in fact its singular locus is $SL(3)$ invariant and thus
it has to be $\iso \P^2$. The normal bundle of this $\P^2$ in $X$ has to be homogeneous, 
as noticed in the lemma (5.2). This cannot occur because there is a 
description of the possible non normal del Pezzo exceptional divisor by Fujita
and the normal bundle of the singular locus (which is $\P^2$) in $X$ is not homogeneous
(see [Fu]).
Alternatively we can use a recent result where, 
together with T. Peternell, we showed that actually there is no non normal
del Pezzo variety of dimension 3 that can be the exceptional locus
of a birational Fano-Mori contraction of a 4-fold
(forcoming paper).

Thus, being smooth, $E$ has to be in the classification of the previous section:
that is $E$ can be either $\P^3$, either $\P(\O_{P^2}(1) \oplus \O_{P^2})$, with conormal 
bundle $\xi \otimes H$ where $\xi$ is the tautological bundle
and $H$ is the pull back of $\O(1)$ from $\P^2$,
or $\P(TP^2)$ with the conormal bundle $\O(1,1)$, the tensor of 
the two line bundles obatined by pulling back $\O(1)$ from the two projection into
$\P^2$.

The case $E = \P^3$ cannot occur because it has a fixed point. 
In the second case we notice that the section at infinity of $\P(\O_{P^2}(1) \oplus \O_{P^2})$
is an orbit $\iso \P^2$ with conormal bundle $N^* =\O(1) \oplus \O(1)$.
Then we can $G$-equivariantly blow-up this orbit and contract the exceptional divisor
into a compact (non projective) manifold which will then 
contains a 1-dimensional orbit, namely
the image of the exceptional divisor isomorphich to $\P^1$; 
this is a contradiction since $SL(3)$ has no
1-dimensional homogeneous variety
(see also the next point c) concerning small contractions). 

The case $E=\P(TP^2)$ can actually occur.
We apply the cone's proposition, thus $X$ is a completion of the
open orbit of the cone $V(\P(TP^2), \O(1,1))$. 
We complete the cone in the trivial way, i.e. 
$X' = \P(\O(1,1) \oplus \O)$. One can easily check that the action
of $SL(3)$ on $X'$ has an open orbit and two closed orbit
isomorphich to $\P(TP^2)$, i.e. the sections at zero and at infinity.
For any other completions $X$ we must have 
a birational map $f : X ---> X'$ induced by identifying the open orbit.
By the theorem of Hironaka we can resolve the indeterminacy of $f$ by blow-ups
along smooth centers. Since the indeterminacy locus is $SL(3)$-equivariant
the centers will be $\P^2$. Let $\sigma : \tilde X \ra X$ be the composition of this blow-ups
and let $\tilde f = f^. \sigma: \tilde X \ra X'$. Since the indeterminacy
of $\tilde f ^{-1}$ is $SL(3)$ equivariant and since $X'$ has no equivariant subset of 
codimension $\geq 2$ then $\tilde f$ is an isomorphism. Thus $X$ is equal
to $X'$ or its smooth $SL(3)$ equivariant blow-down.
One can finally check that in fact there is only one $SL(3)$ equivariant smooth blow-down
of $X'$: namely the zero section is contracted by $\f$ to a point but this contraction
is not elementary and it factors through a smooth blow down with center $\P^2$
(and then through a flop of this $\P^2$ to a point). In particular $X = X'$.

\smallskip \noindent
b) $\f$ is birational, $dim E = 3$ and $dim(\f(E)) > 0$.
In this case, by the $SL(3)$ equivariance of $\f$, we have that
$(\f(E)) = \P^{2}$ and all non trivial fiber are one dimensional.
We can thus apply a result of T. Ando ([An], see also [A-W1]) which says
that in this hypothesis the extremal contraction  
$\f$ is an equivariant smooth blow-up of an orbit $\iso \P^2$
in a smooth manifold $Z$.

\smallskip \noindent
c) $\f$ is a small contractions, i.e. $codim (E) \geq 2$. Thus $E$ has to be of dimension
$2$ and isomorphic to $\P^2$ and with normal bundle 
$N^*$  homogeneus. It is immediate then to check
that $N^* =\O(1) \oplus \O(1)$ (since $det N^* = 2$ and $N^*$ has to be ample); 
this follows also by a general theorem of Kawamata
which describes all small contractions on a smooth 4-fold (see [Ka]).
We blow-up the orbit $\P^2$ and we obtain
a smooth variety with a G-action; but since $N^*$ is ample we can blow-down the 
exceptional divisor in the other direction, i.e. consider the map supported by
$-\hat E - \tau L$ where $L$ is a $\f$-ample divisor and $\tau$ is a rational number
such that $-\hat E - \tau L$ is nef but not ample (thus we can flip the contraction). 
We thus obtain a (smooth) projective variety with a $G$-action and a orbit of dimension one,
the image of the exceptional divisor $\iso \P^1$, a contradiction.  

\medskip
The above three steps prove the following

\proclaim Proposition (5.3). Let $X$ be a smooth projective $4$-fold which has an action
of $SL(3)$ with an open orbit. If $\f :X \ra Z$ is a birational elementary Fano Mori 
contraction then $Z$ is smooth and $\f$ is the blow-up
of an orbit isomorphich to $\P^2$ in $Z$.

In fact in the case in which $E = \P(TP^2)$ we have seen that $\f$ is
not an elementary contractions and it factors through a smooth blow-up (the flop
then is no more a Fano-Mori contraction).

\remark (5.3.1) The above proposition implies that we can run the Minimal Model Program 
{\sl within the category of smooth varieties}. This is true also for the case of quasihomogeneous
3-folds under the action of $SL(2)$ (see [M-U]) and we conjecture it should be true for 
quasihomogeneous $n$-folds under the action of $SL(n-1)$.

\medskip 
Therefore we consider now the cases in which $\f$ is of fiber type.

\smallskip\noindent
d) $\f$ is a conic bundle.

There can be some isolated two dimensional fibers:
then they have to be orbits isomorphic to $\P^2$ and with homogeneous normal bundle.
By the results in [A-W] (in particular 5.9.6) there is only one possibility for
the conormal bundle, namely $N^* = T\P^2(-1)$. 
Moreover in this case $Z$ is smooth thus we use the classification
in the previous section which gives that $Z = \P^3$ (since the images of the isolated exceptional
fibers are fixed points in $Z$).
It is straightforward to see at this point that  $X = \P(T\P^2(-1)\oplus \O)$,
for instance using the results in [B-W].

With the above exception, we have thus that all fibers of
the conic bundle $\f$ are one dimensional; then this implies that $Z$ is smooth,
again by the results in [An],
and we can use the classification in the previous section. 
$Z$ cannot be $\P^3$ since otherwise we will have a one dimensional orbit
(the fiber over the fixed point). Thus $Z= \P( \O_{\P^2}(m) \oplus \O)$ or $\P(TP^2)$;
the first cannot happen since in this case we will not have a dense orbit
while in the second case $X_{(p,q)}= \P(L_{p,q} \oplus \O)$ where $L_{p,q}$ is the line
bundle which corresponds to the character defined on $B$, the Borel subgroup
of $SL(3)$, by

$$\left[\matrix{
a& *&*\cr 
0& e&*\cr 
0&0&i }\right] \ra a^p e^q$$
(for further details see [Na3]).

\smallskip \noindent
e) $\f$ is a Fano fibration over $\P^{2}$; thus it is 
actually an equivariant $\P^2$-bundle, i.e. $X = \P(\E)$ with $\E$ an homogeneous bundle of rank 3
on $\P^2$.
The homogeneous bundles $\O(a) \oplus \O(b) \oplus \O$ don't give a  
quasi homogeneous variety, i.e. there is no open orbit,
except if $a=b= 0$ in which case we have the diagonal action
on $\P^2 \times \P^2$ which has an open orbit.
Therefore $X$ is one of the manifolds $Y_a := \P(T_{\P^2}(a)) \oplus \O)$ or
$\P(S^2T_{\P^2})$.

\bigskip \noindent
To conclude the programm one should now classify all
$X$ Fano fourfolds with $Pic(X) = \Z$
which are quasihomogeneous with the respect to $SL(3)$.

It turns out that the only possibility is the 4-dimensional quadric
embedded into $\P^5= \P(\C^6)$ with the action of $SL(3)$ on $\C^3 \times \C^3$ 
restricted to the quadric.
We are going to present here just a sketch of an argument to prove this fact;
of course it can be checked on the classification of Nakano.

The idea of a direct proof should be the following: 
if $\Omega$ denotes the open orbit then $A:= X \setminus \Omega$
is not empty (see 5.1). One can then prove that $A$ is not connected;
our argument at the moment is long and intricated and we are working to have a clean one
which should be valid in the higher dimensional case too (we expect to have no Fano n-folds 
with $Pic(X) = \Z$ which are quasihomogeneous with the respect to $SL(n-1)$ if $n \geq 5$). 

If $A$ is disconnected we can apply the theorem 2 of Ahiezer (see [Ah]). In particular 
we obtain that $A$ has two components which are homogeneous and isomorphich to $\P^2$ 
with normal bundle $T\P^2(-1)$. Moreover, in his notation, there exists
an equivariant surjective morphism $f: M(SL(3), P, \f) \ra X$, where
in our case $P$ is the parabolic subgroup of $SL(3)$ corresponding
to the first node of the Dynkin diagram and $\f$ is the irreducible
representation corresponding to the tangent bundle twisted by $\O(-1)$. 
It follows now easily that $X$ is isomorphic to the 4-dimensional quadric 
(with the above described $SL(3)$-action) 
and that the map $f$ is the blow-up of one of the two dimensional
orbits. 

Note also that if we blow-up one of the two orbits isomorphic to 
$\P^2$ in the quadric we get a variety isomorphich to $Y_{(-1)}$; if we then blow-up
the strict transform of the other orbit we get a variety isomorphich to
$X_{(0,1)}$.   
  
\bigskip
Collecting all the above together we have thus reproved the following 
\proclaim Theorem (5.4). ([Na3])
Let $X$ be a smooth 4-fold on which $G =SL(3)$ acts with an open orbit.
Then $X$ is isomorphic to one of the following:
\item{1)} $X_{(p,q)} = \P(L_{p,q} \oplus \O)$ where $L_{p,q}$ is the line
bundle on $Z = SL(3) /B = \P (T\P^2)$ which corresponds to the character defined on $B$ by
$\left[\matrix{
a& *&*\cr 
0& e&*\cr 
0&0&i }\right]   \ra a^p e^q.$
\item{2)} $Y_{(a)} = \P(T_{\P^2}(a)) \oplus \O)$ 
\item{3)} $ X = \P(S^2T_{\P^2})$
\item{4)} $X = \P^2 \times \P^2$ and $Bl _{\Delta} (\P^2 \times \P^2)$
\item{5)} $X = \Q^4 \subset \P^5$.

\remark (5.5). The only possible birational elementary contractions between two fourfolds
as in the theorem are blow-ups along smooth orbit isomorphich to $\P^2$, as noticed in (5.3). 
They can now easily be described; the interested reader can look at the proposition 7 of [Na3].

\beginsection References.

\smalltype{

\item{[Ah]} Ahiezer D.N. Dense orbits with two ends, Izv. Akad. Nauk. SSSR
Ser Mat {\bf 41} (1977), 308-324., English Transl Math USSR-Izv {\bf 11} (1997), 293-307.

\item{[An]} Ando, T., On extremal rays of the higher
dimensional varieties, Invent. Math. {\bf 81} (1985), 347---357.

\item {[A-W1]} Andreatta, M, Wi\'sniewski, J.A., A note on non vanishing
and applications, Duke Math. J., {\bf 72} (1993), 739-755.

\item{[A-W2]} Andreatta, M., Wi\'sniewski, J.A., On contractions 
of smooth varieties, Journal of Algebraic Geometry,{\bf 7} (1998), 253-312. 

\item{[B-W]} Ballico, E., Wi\'sniewski, J.A., On B\v anic\v a sheaves and
Fano manifolds, Compositio Math., {\bf 102} (1996), 313-335.

\item{[Bo]} Borel, A., Linear Algebraic Groups, Benjamin (1969).

\item{[Fu]} Fujita, T., On singular del Pezzo varieties, Proc. Conf. 
on Algebraic Geometry L'Aquila 1988, Spriger Verlag {\bf 1417} 1990,117-128.

\item{[F-H]} Fulton, W., Harris, J., Representation theory, 
Graduate texts in mathematics {\bf 129} (1991), Springer Verlag.

\item{[H-O]} Huckleberry A., Oeljeklaus A. A characterization of Complex Homogeneous Cones, 
Math. Zeit. {\bf 170} (1980), 181-194.

\item{[Ka2]} Kawamata, Y., Small contractions of four dimensional
algebraic manifolds, Math. Ann., {\bf 284} (1989), 595-600.

\item{[Ke]} Kebekus, S., Simple models of quasihomogeneous projective
3-folds, Relatively minimal quasihomogeneous projective 3-folds,
On the classification of 3-dimensional $SL(2)$-varieties, 
Alg.Geom. e-prints (1998).

\item{[Ko]} Koll\'ar, J., Rational Curves on Algebraic Varieties, Ergebnisse der Math.,
{\bf 32} (1996), Springer.

\item{[Lu]} Luna, D., Slices \'etales, Bull. Soc. Math. France, Memoire {\bf 33}
(1973), 81-105.

\item{[Ma1-2-3]} Mabuchi,T.,
Equivariant embeddings of normal bundles of fixed point loci,
On the classification of essentially effective $SL(2) \times SL(2)$-actions of algebraic threefolds,
On the classification of essentially effective $SL(n)$-actions on algebraic n-folds,
 Osaka J. Math. {\bf 16} (1979),707-757.

\item{[M-U]} Mukai, S., Umemura H., Minimal rational threefolds, Lect. Notes in Math.
{\bf 1016}, Springer Verlag (1983), 490-518.

\item{[Na1]} Nakano, T.,Regular actions of simple algebraic groups on projective
threefolds, Nagoya Math. J., {\bf 116} (1989), 139-148.

\item{[Na2]} Nakano, T., On equivariant completions of 3-dimensional
homogeneus spaces of $SL(2,\C)$, Japanese J. of Math, {\bf 15} (1989), 221-273.

\item{[Na3]} Nakano, T., On quasi-homogeneus fourfolds of $SL(3,\C)$, 
Osaka J. Math. {\bf 29} (1992),719-733.

}

\end